\documentclass[twoside,11pt,leqno]{article}
\usepackage{amssymb}
\usepackage{amsthm}
\usepackage[tbtags]{amsmath}
\usepackage{doc}
\usepackage{latexsym}
\usepackage{eucal}
\usepackage{enumerate}

\theoremstyle{plain}
\newtheorem*{thm A}{Theorem~A}
\newtheorem*{thm B}{Theorem~C}
\newtheorem*{thm C}{Theorem~D}
\newtheorem*{Main Theorem}{Main Theorem}
\newtheorem*{pro A}{Proposition~E}
\newtheorem*{pro B}{Proposition~F}
\newtheorem*{proposition A}{Theorem~B}

\numberwithin{equation}{section}

\theoremstyle{definition}

\newcommand{\psum}{{+_{\negthinspace\kern-2pt p}}\,}

\newcommand{\R}{\mathbb{R}}
\newcommand{\C}{\mathbb{C}}
\newcommand{\HH}{\mathbb{H}}

\newcommand{\K}{\mathbb{K}}

\newcommand{\id}{\mathrm{id}}

\newcommand{\SO}{\mathrm{SO}}
\newcommand{\End}{\mathrm{End}}

\newcommand{\rk}{\mathrm{rk}}

\newcommand{\eps}{\varepsilon}
\newcommand{\vi}{\varphi}

\newcommand{\Ug}{\mathrm{U}}
\newcommand{\SU}{\mathrm{SU}}
\newcommand{\Sp}{\mathrm{Sp}}

\newcommand{\RE}{\mathop{\mathrm{Re}}\nolimits}
\newcommand{\IM}{\mathop{\mathrm{Im}}\nolimits}

\newcommand{\tr}{\mathop{\mathrm{tr}}\nolimits}

\newcommand{\lieso}{\mathfrak{so}}

\newcommand{\CP}{\ensuremath{\C\mathrm{P}}}
\newcommand{\HP}{\ensuremath{\HH\mathrm{P}}}

\begin{document}
\small{\addtocounter{page}{0} \pagestyle{plain}
\noindent{\scriptsize Proceedings of The 21st International Workshop on \\
Hermitian Symmetric Spaces and Submanifolds 21(2017) 1-2 }

\vspace{0.35in}

\renewcommand{\baselinestretch}{1.5}
\renewcommand{\theequation}{\thesection.\arabic{equation}}
 \makeatletter
 \renewcommand{\section}{\@startsection {section}{0}{0mm}
                        {\baselineskip}{0.3\baselineskip}{\normalfont\bfseries}}
                        \makeatother
                       \makeatother
\makeatletter
   \renewcommand{\subsection}{\@startsection {subsection}{0}{0mm}
                        {\baselineskip}{0.3\baselineskip}{\normalfont\bfseries}}
                        \makeatother

\noindent{\large\bf Low-dimensional totally geodesic submanifolds in ``skew'' position in the symmetric spaces of rank 2}
\vspace{0.15in}\\
\noindent{\sc Sebastian Klein }
\newline
{\it Institut f\"ur Mathematik, Seminargeb\"aude A5, Universit\"at Mannheim, 68131 Mannheim, Germany} \\
{\it e-mail} : { \verb|s.klein@math.uni-mannheim.de| }
\vspace{0.15in}
%
\\
\noindent{\footnotesize(2010 Mathematics Subject Classification : 53C35, 53C40, 53C55, 22E46.)}

\footnote{
\\
* This work was supported by grant Proj. No. ************** from National Research Foundation of Korea.}

\vspace{0.15in}

{\footnotesize
  \noindent Abstract. We use the Cartan representations of \,$\SO(3)$\, and \,$\SU(3)$\,, and an irreducible 14-dimensional representation of \,$\Sp(3)$\, to construct certain
totally geodesic submanifolds in ``skew'' position in the complex quadric, the complex 2-Grassmannians, and the quaternionic 2-Grassmannians.}

\vspace{0.3in}

\section {Introduction}
\setcounter{equation}{0}
\renewcommand{\theequation}{1.\arabic{equation}}
\label{Se:intro}
\vspace{2mm}

There are three infinite series of irreducible, simply connected rank~2 Riemannian symmetric spaces of compact type, namely the complex quadrics \,$Q^m = \SO(m+2)/\SO(2) \times \SO(m)$\,
(which are isomorphic to the real Grassmannians of oriented, real 2-planes \,$G_2^+(\R^{m+2})$\,), the Grassmannians of complex 2-planes \,$G_2(\C^{m+2}) = \SU(m+2) / \mathrm{S}(\Ug(2)\times \Ug(m))$\, and the Grassmannians
of quaternionic 2-planes \,$G_2(\HH^{m+2}) = \Sp(m+2) / \Sp(2)\times\Sp(m)$\,. \,$Q^m$\, is a Hermitian symmetric space. \,$G_2(\C^{m+2})$\, has both a Hermitian and a quaternion Hermitian structure,
and the two are compatible with each other. \,$G_2(\HH^{m+2})$\, admits neither a Hermitian nor a quaternion Hermitian structure. 
Additionally to these three series, there is a finite number of other rank~2 Riemannian symmetric spaces, which we will not consider here. 

When one studies the geometry of a Riemannian symmetric space, one point of significant interest is the study of its totally geodesic submanifolds. In the rank~2 Riemannian symmetric spaces of compact type,
the totally geodesic submanifolds have been completely classified by the author in \cite{Klein1}, \cite{Klein2}, \cite{Klein3}. Most of the totally geodesic submanifolds found in the classification are those that
one would expect from the structure of the root system of the respective symmetric space, and it is not too difficult to write down totally geodesic, isometric embeddings for them explicitly,
see \cite[Section~2]{Klein4}.

However it turns out that the 2-Grassmannians are also contain certain totally geodesic submanifolds of rank 1 of low dimension in a ``skew'' position, which are more unexpected in the sense
that their existence is, perhaps, less obvious from the root diagram of the ambient space. They are
\begin{itemize}
\item a 2-sphere \,$S^2$\, in \,$Q^m$\, for \,$m\geq 3$\,, which is maximal for \,$m=3$\,, denoted by type (A) in \cite[Theorem~4.1]{Klein1};
\item a complex projective space \,$\CP^2$\, in \,$G_2(\C^{m+2})$\, for \,$m\geq 4$\,, which is maximal for \,$m=4$\,, denoted by type \,$(\mathbb{P}, \vi=\arctan(\tfrac12), (\C,2))$\,
  in \cite[Theorem~7.1]{Klein2};
\item a quaternion projective space \,$\HP^2$\, in \,$G_2(\HH^{m+2})$\, for \,$m\geq 5$\,, which is maximal for \,$m=5$\,, denoted by type \,$(\mathbb{P}, \vi=\arctan(\tfrac12), (\HH,2))$\,
  in \cite[Theorem~5.3]{Klein2}.
\end{itemize}
They have a ``skew'' position in their ambient space in the sense that the 2-sphere in the Hermitian symmetric space \,$Q^m$\, is neither complex nor totally real, and the \,$\CP^2$\, in
\,$G_2(\C^{m+2})$\, is neither complex nor totally real with respect to the Hermitian structure, and it is neither quaternionic, nor totally complex nor totally real with respect to the quaternionic
Hermitian structure. The ``characteristic angle'' of the tangent vectors to these totally geodesic submanifolds (referring to the position of these vectors in the Weyl chamber of the
ambient rank 2 symmetric space, see \cite[Section~3]{Reckziegel}) equals \,$\vi=\arctan(\tfrac12)$\,, a value that does not occur for any other totally geodesic submanifolds.
Unlike for other totally geodesic submanifolds of the rank 2 symmetric spaces, in these cases it is not immediately clear 
from the description of the tangent spaces of the totally geodesic submanifolds (in \cite[Theorem~4.1, Type~(A)]{Klein1} and \cite[Theorem~5.3, Type~$(\mathbb{P}, \vi=\arctan(\tfrac12), (\K,2))$]{Klein2})
how a corresponding totally geodesic, isometric embedding can be constructed in a geometrically meaningful way. It is the purpose of the present note to describe such embeddings,
and thereby to explicate the extrinsic geometry of these ``skew'' totally geodesic submanifolds. 

\paragraph{Acknowledgement.} The author is indebted to Professor J.-H.~Eschenburg (Universit\"at Augsburg) and to Professor J.~Berndt (King's College, London) for fruitful discussions
on the matters covered in this note. 

\vspace{2mm}

\section{A totally geodesic 2-sphere in the complex quadric}
\label{Se:Q3}
\setcounter{equation}{0}
\renewcommand{\theequation}{2.\arabic{equation}}
\vspace{2mm}

Let \,$V$\, be a complex linear space of dimension \,$m+2$\,, where \,$m\geq 1$\,, equipped with a positive definite Hermitian product, which we denote by \,$H( \,\cdot\,,\,\cdot\,)$\,.
For \,$z \in V \setminus \{0\}$\, we denote by \,$[z] := \C z$\, the complex line through \,$z$\, and the origin. \,$\CP(V) := \bigr\{\,[z]\,|\,z\in V\setminus \{0\} \,\bigr\}$\, is called
the \emph{complex projective space} over \,$V$\,; it can be made into a \,$(m+1)$-dimensional Hermitian symmetric space of rank \,$1$\, isomorphic to \,$\CP^{m+1}=\SU(m+2)/\mathrm{S}(\Ug(1)\times \Ug(m+1))$\,
in the usual way. If we carry out this construction
with different Hermitian products \,$H$\,, we obtain isometric complex projective spaces. 

Now suppose that a non-degenerate, symmetric, complex-bilinear form \,$\beta: V \times V \to \C$\, is given. We say that \,$\beta$\, is \emph{compatible} with the Hermitian product \,$H$\,, if
the anti-linear map \,$A: V \to V$\, determined by the equation \,$H(v,Aw)=\beta(v,w)$\, for \,$v,w\in V$\, is a real structure on \,$V$\,, i.e.~if \,$A$\, is an involutive, anti-holomorphic
isometry of \,$(V,H)$\,. If \,$\beta$\, is compatible with \,$H$\,, we call the complex hypersurface
$$ Q(V,\beta) := \bigr\{ \,[z]\,|\, z \in V\setminus\{0\}, \beta(z,z)=0 \,\bigr\} $$
of \,$\CP(V)$\, the \,$m$-dimensional \emph{complex quadric} in \,$\CP(V)$\,. \,$Q(V,\beta)$\, is a Hermitian symmetric space of rank \,$2$\, isomorphic to \,$Q^m = \SO(m+2)/\SO(2) \times \SO(m)$\,.
If we carry out this construction with different
symmetric bilinear forms \,$\beta$\, which are compatible with \,$H$\,, we obtain isometric complex quadrics.

We want to construct a certain 2-sphere \,$M\cong \SO(3)/\SO(2)$\,, which is a maximal totally geodesic submanifold in \,$Q^3 = \SO(5)/\SO(2)\times \SO(3)$\, of type (A)
as in \cite[Theorem~4.1]{Klein1}. Because complete totally geodesic
submanifolds of symmetric spaces are symmetric subspaces, the transvection group of \,$M$\, corresponds to a subgroup of the transvection group \,$\SO(5)$\, of \,$Q^3$\, which is isomorphic to \,$\SO(3)$\,, or
in other words, to a 5-dimensional orthogonal, real representation of \,$\SO(3)$\,. It turns out that the representation which we need is the \emph{Cartan representation} of \,$\SO(3)$\,.

To construct the Cartan representation, consider the 5-dimensional real linear space \,$\End_+^0(\R^3)$\, of self-adjoint, trace-free endomorphisms on \,$\R^3$\,, along with its complexification,
the 5-dimensional complex-linear space \,$V := \End_+^0(\R^3)\otimes \C$\, of symmetric, trace-free complex \,$(3\times 3)$-matrices. The sesqui-linear form
$$ H(X,Y) := \tr(X \cdot \overline{Y}) \quad\text{for}\quad X,Y\in V $$
is a Hermitian product on \,$V$\,, and the symmetric bilinear form
$$ \beta(X,Y) := \tr(X \cdot Y) \quad\text{for}\quad X,Y\in V $$
is non-degenerate and compatible with \,$H$\,. For \,$B \in \SO(3)$\, and \,$X\in V$\,, we have \,$BXB^{-1} \in V$\,, therefore we obtain the 5-dimensional representation of \,$\SO(3)$\,
$$ \SO(3) \times V \to V,\; (B,X) \mapsto BXB^{t}=BXB^{-1} \;, $$
which is called the \emph{Cartan representation}. It leaves the Hermitian product \,$H$\, and the real structure \,$V \to V, \; X \mapsto \overline{X}$\, invariant, therefore it is
an orthogonal, real representation of \,$\SO(3)$\,. Because it moreover leaves the bilinear form \,$\beta$\, invariant, it induces an isometric action of \,$\SO(3)$\, on the 3-dimensional
complex quadric \,$Q(V,\beta)$\,. 

Let
$$ Z_0 := \tfrac12 \left( \begin{smallmatrix} 0 & 1 & i \\ 1 & 0 & 0 \\ i & 0 & 0 \end{smallmatrix} \right) \in V \; . $$
Then we have \,$H(Z_0,Z_0)=1$\, and \,$\beta(Z_0,Z_0)=0$\,, and therefore \,$[Z_0] \in Q(V,\beta)$\,. Let \,$M$\, be the orbit through \,$[Z_0]$\, of the \,$\SO(3)$-action on \,$Q(V,\beta)$\, induced
by the Cartan representation of \,$\SO(3)$\,, and let \,$K\subset \SO(3)$\, be the isotropy subgroup of this action at the point \,$[Z_0]$\,. Then \,$M$\, is as a homogeneous space isomorphic
to \,$\SO(3)/K$\,. In the following paragraph we will show that \,$K$\, is isomorphic to \,$\SO(2)$\,, and therefore \,$M$\, is isomorphic to \,$\SO(3)/\SO(2) = S^2$\,. 

\emph{Proof that \,$K\cong \SO(2)$\,.}
  For given \,$B\in \SO(3)$\, we have \,$B\in K$\, if and only if \,$BZ_0B^{-1} = \lambda Z_0$\, for some \,$\lambda\in S^1$\,. If we denote by \,$e_1, e_2, e_3$\, the standard unit vectors of \,$\R^3$\,,
  we have
  $$ Z_0e_1 = \tfrac12(e_2+ie_3)\;,\quad Z_0e_2 = \tfrac12 e_1 \quad\text{and}\quad Z_0e_3 = \tfrac{i}{2}e_1 \; . $$
  Therefore the condition \,$BZ_0B^{-1}e_1=\lambda Z_0e_1$\, implies \,$Be_1 \in \mathrm{span}\{e_1\}$\, and therefore (because \,$B$\, is orthogonal) \,$Be_1 = \eps e_1$\, for some
  \,$\eps \in \{\pm 1\}$\,. Again because \,$B$\, is orthogonal, and \,$\det(B)=1$\, holds, it follows that there exist \,$r,s\in \R$\, with \,$r^2 + s^2 = 1$\, so that
  $$ Be_2 = re_2 + se_3 \quad\text{and}\quad Be_3 = -\eps s e_2 + \eps r e_3 $$
  holds. From the condition \,$BZ_0B^{-1}e_2 = \lambda Z_0 e_2$\, one now obtains by a straightforward calculation \,$\lambda = \eps\,(r-is)$\,, and then another calculation
  from the condition \,$BZ_0B^{-1}e_3 = \lambda Z_0 e_3 = \eps(r-is)Z_0e_3$\, gives \,$\eps=1$\,. Hence
  $$ B = \left( \begin{smallmatrix} 1 & 0 & 0 \\ 0 & r & -s \\ 0 & s & r \end{smallmatrix} \right) $$
  holds. This shows that \,$K\cong \SO(2)$\,.
  \strut\hfill$\Box$
  
It remains to show that \,$M$\, is a totally geodesic submanifold of \,$Q(V,\beta)$\, which has the ``skew'' position of type $(A)$ in \cite[Theorem~4.1]{Klein1} as described in the Introduction.
To show that \,$M$\, is totally geodesic in \,$Q^m$\,, we consider the Lie group homomorphism induced by the Cartan representation:
$$ \Phi: \SO(3) \to \SO(\End_+^0(\R^3)),\; B \mapsto (X \mapsto BXB^{-1}) \; . $$
We need to show that the linearisation of \,$\Phi$\,
$$ \Phi_*: \lieso(3) \to \lieso(\End_+^0(\R^3)),\; Y \mapsto (X \mapsto YX-XY) $$
respects the Cartan decompositions of the symmetric spaces \,$M$\, and \,$Q(V,\beta)$\,. The Lie algebra of the transvection group of \,$M$\, is \,$\lieso(3)$\, (i.e.~the space of
skew-symmetric, real \,$(3\times 3)$-matrices), and its Cartan decomposition is given by \,$\lieso(3) = \mathfrak{k} \oplus \mathfrak{m}$\,
with
$$ \mathfrak{k} = \{X \in \lieso(3)|Xe_1=0\}\cong \lieso(2) \quad\text{and}\quad \mathfrak{m} = \{X \in \lieso(3)|Xe_2,Xe_3 \in \R\,e_1\} \; . $$
The Lie algebra of the transvection group of \,$Q(V,\beta)$\, is \,$\lieso(\End_+^0(\R^3))$\,, i.e.~the space of skew-adjoint linear maps on the 5-dimensional real linear space \,$\End_+^0(\R^3)$\,
equipped with the inner product given by the real part of \,$H$\,. We consider the splitting \,$\End_+^0(\R^3) = W_1\oplus W_2$\, with
$$ W_1 = \{X \in \End_+^0(\R^3)|Xe_2,Xe_3 \in \R e_1\} \quad\text{and}\quad W_2 = \{X \in \End_+^0(\R^3)|Xe_1 \in \R e_1\} \; ; $$
note that \,$\dim(W_1)=2$\, and \,$\dim(W_2)=3$\, holds. The Cartan decomposition of \,$Q(V,\beta)$\, is then given by \,$\lieso(\End_+^0(\R^3)) = \widetilde{\mathfrak{k}} \oplus \widetilde{\mathfrak{m}}$\, with
\begin{align*}
  \mathfrak{k} & = \bigr\{ \, A \in \lieso(\End_+^0(\R^3)) \,\bigr|\, AW_1 \subset W_1,\, AW_2 \subset W_2 \, \bigr\} \cong \lieso(2) \oplus \lieso(3) \;, \\
  \mathfrak{m} & = \bigr\{ \, A \in \lieso(\End_+^0(\R^3)) \,\bigr|\, AW_1 \subset W_2,\, AW_2 \subset W_1 \, \bigr\} \;.
\end{align*}
Because \,$\Phi$\, induces a well-defined action on \,$Q(V,\beta)$\,, it is clear that \,$\Phi_*\mathfrak{k} \subset \widetilde{\mathfrak{k}}$\, holds, and we need to show that
also \,$\Phi_*\mathfrak{m} \subset \widetilde{\mathfrak{m}}$\, holds. For this purpose, let \,$Y \in \mathfrak{m}$\, be given. Then there exist \,$t,s\in \R$\, so that
$$ Y = \left( \begin{smallmatrix} 0 & t & s \\ -t & 0 & 0 \\ -s & 0 & 0 \end{smallmatrix} \right) $$
holds. An explicit calculation shows that for \,$X \in W_1$\,, say
$$ X = \left( \begin{smallmatrix} 0 & a & b \\ a & 0 & 0 \\ b & 0 & 0 \end{smallmatrix} \right) $$
with \,$a,b\in \R$\,, one has
$$ (\Phi_*Y)X = YX-XY = \left( \begin{smallmatrix} -2ta-2sb & 0 & 0 \\ 0 & 2ta & sa+tb \\ 0 & tb+sa & 2sb \end{smallmatrix} \right) \in W_2 \; . $$
Similarly, for \,$X \in W_2$\,, say
$$ X = \left( \begin{smallmatrix} a & 0 & 0 \\ 0 & b & d \\ 0 & d & c \end{smallmatrix} \right) $$
with \,$a,b,c,d \in \R$\, and \,$a+b+c=0$\,, one has
$$ (\Phi_*Y)X = YX-XY = \left( \begin{smallmatrix} 0 & (a-b)t-ds & -dt + (a-c)s \\ (a-b)t-ds & 0 & 0 \\ -dt + (a-c)s & 0 & 0 \end{smallmatrix} \right) \in W_1 \; . $$
This shows that \,$\Phi_*Y \in \widetilde{\mathfrak{m}}$\, holds. Therefore \,$M$\, is a symmetric subspace, hence a totally geodesic submanifold of \,$Q(V,\beta)$\,.

It follows that for any \,$X\in \mathfrak{m}$\,, the curve \,$\gamma(t) := \exp(tX)\,Z_0\,\exp(tX)^{-1}$\, in \,$\mathbb{S}(V)$\, projects to a geodesic of \,$Q(V,\beta)$\, starting in \,$[Z_0]$\,
and running entirely in \,$M$\,, and every such geodesic is obtained in this way. We have \,$\dot{\gamma}(0) = X\,Z_0 - Z_0\,X$\,, and therefore \,$T_{[Z_0]}M = \{X\,Z_0-Z_0\,X|X\in\mathfrak{m}\}$\,.
Because \,$\mathfrak{m}$\, is spanned by
$$ X_1 := \left( \begin{smallmatrix} 0 & 1 & 0 \\ -1 & 0 & 0 \\ 0 & 0 & 0 \end{smallmatrix} \right) \quad\text{and}\quad  X_2 := \left( \begin{smallmatrix} 0 & 0 & 1 \\ 0 & 0 & 0 \\ -1 & 0 & 0 \end{smallmatrix} \right) \;, $$
\,$T_{[Z_0]}M$\, is spanned by
$$ Y_1 := X_1\,Z_0 - Z_0\,X_1 = \left( \begin{smallmatrix} 1 & 0 & 0 \\ 0 & -1 & -i/2 \\ 0 & -i/2 & 0 \end{smallmatrix} \right) \quad\text{and}\quad
Y_2 := X_2\,Z_0 - Z_0\,X_2 = \left( \begin{smallmatrix} i & 0 & 0 \\ 0 & 0 & -1/2 \\ 0 & -1/2 & -i \end{smallmatrix} \right) \; . $$
Note that we have \,$\RE(Y_k) \perp \IM(Y_k)$\,, 
\,$\|\RE(Y_1)\| = 2\cdot \|\IM(Y_1)\|$\, and \,$\|\IM(Y_2)\| = 2\cdot \|\RE(Y_2)\|$\,. These equations show that both vectors \,$Y_1$\, and \,$Y_2$\, have the ``characteristic
angle'' \,$\vi=\arctan(\tfrac12)$\, in the sense of \cite[Section~3]{Reckziegel} (although they lie in different Weyl chambers). Therefore the totally geodesic submanifold \,$M$\, of
\,$Q(V,\beta)$\, is of type (A) as in \cite[Theorem~4.1]{Klein1}.

\vspace{2mm}

\section{A totally geodesic \,$\CP^2$\, in the complex 2-Grassmannian}
\label{Se:G2C6}
\setcounter{equation}{0}
\renewcommand{\theequation}{3.\arabic{equation}}
\vspace{2mm}

Again let \,$V$\, be a complex linear space of dimension \,$m+2$\,, where \,$m\geq 2$\,, equipped with a positive definite Hermitian product \,$H(\,\cdot\,,\,\cdot\,)$\,.
We denote the set of complex planes (2-dimensional complex linear subspaces of \,$V$\,) by \,$G_2(V)$\,. \,$G_2(V)$\, is called the \emph{complex 2-Grassmannian over \,$V$\,}.
It is well-known that the Hermitian product \,$H$\, induces the structure of a Riemannian manifold on \,$G_2(V)$\,; in this way, \,$G_2(V)$\, becomes a Hermitian symmetric space
that also carries the structure of a quaternionic Hermitian symmetric space compatible with the other structures. It is of complex dimension \,$2m$\,
(real dimension \,$4m$\,), and as a symmetric space of compact type and rank 2,
isomorphic to \,$\SU(m+2)/\mathrm{S}(\Ug(2)\times \Ug(m))$\,.

We want to construct the totally geodesic \,$\CP^2$\, in \,$G_2(\C^6)$\, which is of type \,$(\mathbb{P}, \vi=\arctan(\tfrac12), (\C,2))$\, as in \cite[Theorem~7.1]{Klein2}. Because
\,$\CP^2$\, is isomorphic to \,$\SU(3)/\mathrm{S}(\Ug(1)\times \Ug(2))$\, and \,$G_2(\C^6)$\, is isomorphic to \,$\SU(6)/\mathrm{S}(\Ug(2)\times \Ug(4))$\,, we need a suitable unitary
representation of \,$\SU(3)$\, on a complex-6-dimensional complex linear space \,$V$\, that descends to an action on \,$G_2(V)$\,. This purpose is again filled by the
Cartan representation of \,$\SU(3)$\,.

Let \,$V:=\End_+(\R^3) \otimes \C$\, be the 6-dimensional complex linear space of symmetric (\emph{not} self-adjoint!) \,$(3\times 3)$-matrices, which we again equip with the non-degenerate Hermitian product
$$ H(X,Y) := \tr(X \cdot \overline{Y}) \quad\text{for}\quad X,Y\in V \; . $$
For \,$B\in \SU(3)$\, and \,$X\in V$\, we have \,$BXB^t \in V$\,, and thus we again have the Cartan representation
$$ \SU(3)\times V \to V,\; (B,X) \mapsto BXB^t \; . $$
This representation is unitary with respect to \,$H$\,; indeed, for \,$B\in \SU(3)$\, and \,$X,Y\in V$\, we have
\begin{align*}
  H(BXB^t,BYB^t) & = \tr(BXB^t \cdot \overline{BYB^t}) = \tr(B\,X\,B^t\,(B^t)^{-1}\,\overline{Y}\,B^{-1}) \\
  & = \tr(B\,X\,\overline{Y}\,B^{-1}) = \tr(X\,\overline{Y}) = H(X,Y) \; .
\end{align*}
Note that unlike the Cartan representation of \,$\SO(3)$\, we used in Section~\ref{Se:Q3}, the Cartan representation of \,$\SU(3)$\, does not leave the trace of members of \,$V$\, invariant;
this is the reason why we cannot restrict to trace-free matrices here. The Cartan representation of \,$\SU(3)$\, induces an isometric action of \,$\SU(3)$\, on the 2-Grassmannian
\,$G_2(V)$\,.

Let \,$Z_0 \in G_2(V)$\, be the 2-dimensional complex linear subspace of \,$V$\, spanned by
$$ A_1 := \left( \begin{smallmatrix} 0 & 1 & 0 \\ 1 & 0 & 0 \\ 0 & 0 & 0 \end{smallmatrix} \right) \quad\text{and}\quad A_2 := \left( \begin{smallmatrix} 0 & 0 & 1 \\ 0 & 0 & 0 \\ 1 & 0 & 0 \end{smallmatrix} \right) \; . $$
Let \,$M$\, be the orbit through \,$Z_0$\, of the isometric action of \,$\SU(3)$\, on \,$G_2(V)$\, induced by the Cartan representation, and let
\,$K\subset \SU(3)$\, be the isotropy subgroup of this action at the point \,$Z_0$\,. Then \,$M$\, is as a homogeneous space isomorphic to \,$\SU(3)/K$\,. We will show below
that \,$K$\, is isomorphic to \,$\mathrm{S}(\Ug(1)\times \Ug(2))$\,, and therefore \,$M$\, is isomorphic to \,$\SU(3)/\mathrm{S}(\Ug(1)\times \Ug(2)) \cong \CP^2$\,.
Analogous calculations as in Section~\ref{Se:Q3} show that the Lie group homomorphism
$$ \Phi: \SU(3) \to \SU(V),\; B \mapsto (X \mapsto BXB^t) $$
induces a homomorphism of symmetric spaces, and therefore \,$M$\, is a symmetric subspace, hence a totally geodesic submanifold of \,$G_2(V)$\,. Moreover, one can calculate
the tangent space of \,$M$\, like in Section~\ref{Se:Q3}, and from this calculation it follows that \,$T_pM$\, has a 2-dimensional totally real subspace consisting entirely
of vectors of ``characteristic angle'' \,$\vi=\arctan(\tfrac12)$\,. Therefore the totally geodesic submanifold \,$M$\, of \,$G_2(V)$\, is of type
\,$(\mathbb{P}, \vi=\arctan(\tfrac12), (\C,2))$\, as  in \cite[Theorem~7.1]{Klein2}.

\emph{Proof that \,$K\cong \mathrm{S}(\Ug(1)\times \Ug(2))$\,.}
Let \,$B \in \SU(3)$\, be given. \,$(A_1,A_2)$\, is a unitary basis of \,$Z_0$\,, and we have \,$B\in K$\, if and only if \,$(BA_1B^t,BA_2B^t)$\, is another unitary basis of \,$Z_0$\,.
The latter condition means that there exist \,$a,b\in \C$\, with \,$|a|^2+|b|^2=1$\, so that
$$ BA_1B^t = a\cdot A_1 + b\cdot A_2 \quad\text{and}\quad BA_2B^t = -\overline{b}\cdot A_1 + \overline{a}\cdot A_2 $$
holds. Let us write \,$B=(b_{k\ell})_{k,\ell=1,2,3}$\,. Then evaluating the preceding equations and comparing matrix entries
gives that the following equations hold:
\begin{gather}
  \label{eq:G2Cn:eq1}
  b_{11}\,b_{12} = b_{11}\,b_{13} = b_{21}\,b_{22}=b_{21}\,b_{23} = b_{31}\,b_{32}=b_{31}\,b_{33} = 0 \\
  \label{eq:G2Cn:eq2}  
  b_{22}\,b_{31} + b_{21}\,b_{32} = b_{23}\,b_{31} + b_{21}\,b_{33} = 0 \; .
\end{gather}
Equations~\eqref{eq:G2Cn:eq1} show that either \,$b_{21}=0$\, or \,$b_{22}=b_{23}=0$\, holds. If \,$b_{22}=b_{23}=0$\, were true, we would have \,$b_{21}\neq 0$\, (because of \,$\rk(B)=3$\,), and therefore
Equations~\eqref{eq:G2Cn:eq2} would imply \,$b_{32}=b_{33}=0$\,. Therefore we would have \,$Be_2,Be_3 \in \C e_1$\,, which is a contradiction to \,$\rk(B)=3$\,. Therefore we have
\,$b_{21}=0$\,. Similarly, \,$b_{31}=0$\, also holds. Because of \,$\rk(B)=3$\, we then have \,$b_{11}\neq 0$\,. Equations~\eqref{eq:G2Cn:eq1} then show \,$b_{12}=b_{13}=0$\,. Therefore \,$B$\,
leaves \,$\C e_1$\, and \,$\C e_2 \oplus \C e_3$\, invariant, and thus we have \,$B \in \mathrm{S}(\Ug(1) \times \Ug(2))$\,. 
\strut\hfill$\Box$

\vspace{2mm}
\section{A totally geodesic \,$\HP^2$\, in the quaternionic 2-Grassmannian}
\label{Se:G2H7}
\setcounter{equation}{0}
\renewcommand{\theequation}{4.\arabic{equation}}
\vspace{2mm}

Let \,$V$\, be a \,$2(m+2)$-dimensional complex vector space equipped with a Hermitian product \,$H(\,\cdot\,,\,\cdot\,)$\,. A \emph{quaternionic structure} on \,$V$\, is an
anti-linear isometry \,$J:V \to V$\, with \,$J^2=-\id$\,, and we call \,$V$\, a \,$(m+2)$-dimensional \emph{quaternionic linear space} if it equipped with a quaterionic structure \,$J$\,.
By identifying the application of \,$J$\, with the multiplication with an imaginary unit quaternion \,$j$\, orthogonal to \,$i\in \C$\,, \,$V$\, then indeed becomes a
module over the skew-field of quaternions \,$\HH$\,, isomorphic to \,$\HH^{m+2}$\,. In this setting
\begin{equation}
\label{eq:G2H7:omega}
\omega(v,w) := H(v,Jw)
\end{equation}
is a non-degenerate, complex-valued symplectic form on \,$V$\,, and the group
$$ \Sp(V) = \{B\in \Ug(V)|\forall v,w \in V: \omega(Bv,Bw)=\omega(v,w)\} = \{B \in \Ug(V) | B \circ J = J \circ B \} $$
is called the \emph{symplectic group} over \,$V$\,; its members are called the \emph{symplectic maps} of \,$V$\,. We automatically have \,$\Sp(V) \subset \SU(V)$\,.

A \emph{quaternionic plane} in \,$V$\, is a complex-4-dimensional linear subspace of \,$V$\, that is invariant under \,$J$\,. As a \,$\HH$-module, it has the dimension \,$2$\,.
The \emph{quaternionic 2-Grassmannian} \,$G_2(V)$\, over \,$V$\, is the set of all quaternionic planes of \,$V$\,. As before, the Hermitian product \,$H$\, induces the structure of a
Riemannian symmetric space on \,$G_2(V)$\,, it is of compact type and rank \,$2$\,, has real dimension \,$8m$\,, and is as a symmetric space isomorphic to \,$\Sp(m+2)/\Sp(2)\times \Sp(m)$\,.
Note that on \,$G_2(V)$\,, there neither exists a compatible Hermitian structure, nor a compatible quaternionic Hermitian structure.

Here we want to construct the totally geodesic \,$\HP^2$\, in \,$G_2(\HH^7)$\, that is of type \,$(\mathbb{P}, \vi=\arctan(\tfrac12), (\HH,2))$\, as  in \cite[Theorem~5.3]{Klein2}.
Because \,$\HP^2$\, is isomorphic to \,$\Sp(3)/\Sp(1)\times \Sp(2)$\, and \,$G_2(\HH^7)$\, is isomorphic to \,$\Sp(7)/\Sp(2)\times \Sp(5)$\,, we need a suitable
quaternionic representation of \,$\Sp(3)$\, on a complex-14-dimensional linear space \,$V$\,, that descends to an action on the quaternionic 2-Grassmannian over \,$V$\,.

There is exactly one such representation of \,$\Sp(3)$\, (see \cite[Chapter~VI, Section~(5.3), p.~269ff.]{BT}). To construct it, we unfortunately cannot use a Cartan
representation as in the previous two sections. Rather, we note that the vector representation of \,$\Sp(3)$\, on \,$\C^6$\, induces a representation of \,$\Sp(3)$\,
on the three-fold exterior product \,$\bigwedge^3 \C^6$\,. It turns out that this 20-dimensional representation decomposes into two irreducible components \,$V_1$\, and \,$V_2$\,:
One of them, say \,$V_1$\,, is 6-dimensional, and equivalent to the vector representation of \,$\Sp(3)$\, on \,$\C^6$\,; the other, \,$V_2$\,, is therefore 14-dimensional,
and the representation that we seek.

To construct this representation explicitly, we equip \,$\C^6$\, with the quaternionic structure \,$J:\C^6 \to \C^6$\, characterized by
$$ Je_1 = e_4, Je_2=e_5, Je_3 = e_6, Je_4 = -e_1, Je_5=-e_2 \text{ and } Je_6 = -e_3 \;, $$
where \,$e_1,\dotsc,e_6$\, are the usual standard basis vectors of \,$\C^6$\,, and thereby with the associated symplectic form \,$\omega$\, given by Equation~\eqref{eq:G2H7:omega}.
We consider the three-fold exterior product \,$\bigwedge^3 \C^6$\, of \,$\C^6$\,; this is the 20-dimensional
complex vector space with the basis
$$ \{e_j \wedge e_k \wedge e_\ell | 1 \leq j < k < \ell \leq 6\} \;, $$
where the tri-linear map \,$(\C^6)^3 \to \bigwedge^3\C^6, \; (v_1,v_2,v_3) \mapsto v_1\wedge v_2 \wedge v_3$\, is alternating. Note that any \,$\C$-linear (anti-linear) map \,$A: \C^6 \to \C^6$\,
induces a \,$\C$-linear (anti-linear) map \,$A^{(3)}: \bigwedge^3 \C^6 \to \bigwedge^3 \C^6$\, characterized by
$$ \forall v_1,v_2,v_3 \in \C^6 \; : \; A^{(3)}(v_1\wedge v_2 \wedge v_3) = Av_1 \wedge Av_2 \wedge Av_3 \; . $$
\,$\bigwedge^3 \C^6$\, becomes a quaternionic linear space via the Hermitian product \,$\widehat{H}$\, characterized by
$$ \widehat{H}(u_1\wedge u_2\wedge u_3, v_1\wedge v_2 \wedge v_3) = \det(H(u_k,v_\ell))_{k,\ell=1,2,3} \quad\text{for}\quad u_k,v_\ell \in \C^6$$
and the quaternionic structure given by \,$J^{(3)}$\,.
Note that for \,$B\in \Sp(\C^6)=\Sp(3)$\,, we have \,$B^{(3)} \in \Sp(\bigwedge^3\C^6)$\,, and therefore we obtain a quaternionic representation of \,$\Sp(3)$\, on \,$\bigwedge^3\C^6$\, by
$$ \Sp(3) \times {\textstyle \bigwedge^3\C^6} \to {\textstyle \bigwedge^3\C^6},\; (B,\xi) \mapsto B^{(3)}\xi \; . $$
The \emph{contraction} of \,$\bigwedge^3\C^6$\, is the linear map \,$\varkappa: \bigwedge^3\C^6 \to \C^6$\, characterized by
$$  \varkappa(v_1\wedge v_2 \wedge v_3) = \omega(v_1,v_2)\,v_3 + \omega(v_2,v_3)\,v_1 + \omega(v_3,v_1)\,v_2 \quad\text{for}\quad v_1,v_2,v_3 \in \C^6 \; . $$
\,$\varkappa$\, is surjective, \,$J$-invariant (i.e.~\,$\varkappa \circ J^{(3)} = J \circ \varkappa$\,) and \,$\Sp(3)$-invariant; the latter property means that we have
$$ \varkappa(B^{(3)}\xi) = B(\varkappa(\xi)) \quad\text{for}\quad B\in \Sp(3),\; \xi \in {\textstyle \bigwedge^3\C^6} \; . $$
Therefore \,$V_1 := \ker \varkappa$\, and \,$V_2 := V_1^\perp$\, are complementary quaternionic subspaces of \,$\bigwedge^3\C^6$\, with \,$\dim_{\C}(V_1)=14$\,, \,$\dim_{\C}(V_2)=6$\,,
which are both invariant under the \,$\Sp(3)$-action on \,$\bigwedge^3\C^6$\,. The representation of \,$\Sp(3)$\, on \,$V_2$\, is equivalent to the ``vector representation'' of \,$\Sp(3)$\,,
i.e.~to the natural representation of \,$\Sp(3)$\, on \,$\C^6$\,. The representation of \,$\Sp(3)$\, on \,$V_1$\, turns out to be irreducible; it is one of the fundamental representations
of \,$\Sp(3)$\,, see \cite[p.~271f.]{BT}. We will use this representation to construct the desired totally geodesic \,$\HP^2$\, in \,$G_2(V_1)$\,. Note that \,$V_1$\, has
the quaternionic dimension \,$7$\,; as a \,$\HH$-module, it is spanned by the following seven vectors:
\begin{gather*}
  e_1 \wedge e_2 \wedge e_3,\quad Je_1 \wedge e_2 \wedge e_3,\quad e_1 \wedge Je_2 \wedge e_3,\quad e_1 \wedge e_2 \wedge Je_3, \\
  e_1 \wedge e_2 \wedge Je_2 - e_1 \wedge e_3 \wedge Je_3,\\
  e_2 \wedge e_1 \wedge Je_1 - e_2 \wedge e_3 \wedge Je_3, \\
  e_3 \wedge e_1 \wedge Je_1 - e_3 \wedge e_2 \wedge Je_2.
\end{gather*}

Clearly the quaternionic representation of \,$\Sp(3)$\, on \,$V_1$\, which we just described induces an action of \,$\Sp(3)$\, on the quaternionic 2-Grassmannian \,$G_2(V_1)$\, by isometries. Let 
$$ \xi_1 := e_2 \wedge e_1 \wedge Je_1 - e_2 \wedge e_3 \wedge Je_3 \quad\text{and}\quad \xi_2 := e_3 \wedge e_1 \wedge Je_1 - e_3 \wedge e_2 \wedge Je_2   \; . $$
Then we have \,$\varkappa(\xi_1)=\varkappa(\xi_2)=0$\, and therefore \,$Z_0 := \mathrm{span}_{\C}\{\xi_1,J^{(3)}\xi_1,\xi_2,J^{(3)}\xi_2\} \in G_2(V_1)$\,.
We let \,$M$\, be the orbit through \,$Z_0$\, of the \,$\Sp(3)$-action on \,$G_2(V_1)$\,. Then
\,$M$\, is as a homogeneous space isomorphic to \,$\Sp(3)/K$\,, where \,$K \subset \Sp(3)$\, is the isotropy group at \,$Z_0$\, of the \,$\Sp(3)$-action on \,$G_2(V_1)$\,.
One can show that \,$K$\, is isomorphic to \,$\Sp(1)\times \Sp(2)$\,. 
Therefore \,$M$\, is as a homogeneous space isomorphic to \,$\Sp(3)/\Sp(1)\times \Sp(2) \cong \HP^2$\,. Similar calculations as in Section~\ref{Se:Q3}
again show that \,$M$\, is a symmetric subspace, hence a totally geodesic submanifold of \,$G_2(V_1)$\,. Also in a similar way, the tangent space to \,$M$\, at \,$Z_0$\, can be computed,
showing by comparision with \cite[Theorem~5.3]{Klein2} that the totally geodesic submanifold \,$M$\, is of type \,$(\mathbb{P}, \vi=\arctan(\tfrac12), (\HH,2))$\, as in that theorem.

The Cartan representations of \,$\SO(3)$\, resp.~\,$\SU(3)$\, that we used in the constructions of Sections~\ref{Se:Q3} and \ref{Se:G2C6} can be recovered from the
representation of \,$\Sp(3)$\, on \,$V_1$\, by appropriate restriction.
In fact, if we restrict the representation of \,$\Sp(3)$\, on \,$V_1$\, to \,$\SU(3)\subset \Sp(3)$\,,
this representation leaves a totally complex, complex-6-dimensional vector space \,$V_{\C}$\, invariant. The action of \,$\SU(3)$\, on \,$V_{\C}$\, is equivalent to the Cartan
representation of \,$\SU(3)$\, we considered in Section~\ref{Se:G2C6}. Moreover, if we further restrict the representation of \,$\SU(3)$\, on \,$V_{\C}$\, to \,$\SO(3) \subset \SU(3)$\,,
it turns out that we obtain a real representation on \,$V_{\C}$\,, which leaves a complex-1-dimensional subspace of \,$V_{\C}$\, invariant. If we denote the complementary complex-5-dimensional
subspace of \,$V_{\C}$\, by \,$V_{\R}$\,, then the resulting representation of \,$\SO(3)$\, on \,$V_{\R}$\, is equivalent to the Cartan representation of \,$\SO(3)$\, which we used
in Section~\ref{Se:Q3}. 

\vspace{2mm}


}
\end{document}